\documentclass{amsart}

\newtheorem{theorem}{Theorem}
\newcommand{\bt}{\begin{theorem}}
\newcommand{\et}{\end{theorem}}
\newtheorem{lemma}{Lemma}
\newcommand{\bl}{\begin{lemma}}
\newcommand{\el}{\end{lemma}}
\newtheorem{corollary}{Corollary}
\newcommand{\bc}{\begin{corollary}}
\newcommand{\ec}{\end{corollary}}
\newcommand{\beq}{\begin{equation}}
\newcommand{\eeq}{\end{equation}}
\newcommand{\benum}{\begin{enumerate}}
\newcommand{\eenum}{\end{enumerate}}

\begin{document}

\title{Paul Erd\H os and  additive bases}
\author{Melvyn B. Nathanson}
\address{Department of Mathematics\\
Lehman College (CUNY)\\
Bronx, NY 10468}
\email{melvyn.nathanson@lehman.cuny.edu}

\subjclass[2010]{.}
\keywords{Additive basis, thin basis, minimal basis, maximal nonbasis, essential component, erd\H os-Tur\' an conjecture.}

\begin{abstract} This is a survey of some of Erd\H os's work on bases in additive number theory.  \end{abstract}

\maketitle

\section{Additive bases}
Paul Erd\H os, while he was still in his 20s, wrote a series of extraordinarily 
beautiful papers in additive and combinatorial number theory.  
The key concept is  \emph{additive basis}.

Let $A$ be a set of nonnegative integers, let $h$ be a positive integer,
and let $hA$ denote the set of integers 
that can be represented as the sum of exactly $h$ elements of $A$, with repetitions allowed.  
A central problem in additive number theory is to describe the sumset $hA$.  
The set $A$ is called an \emph{additive basis of order $h$} if every nonnegative integer 
can be represented as the sum of exactly $h$ elements of $A$.  
For example, the set of squares is a basis of order 4 (Lagrange's theorem) 
and the set of nonnegative cubes is a basis of order 9 (Wieferich's theorem).  

The set $A$ of nonnegative integers is an \emph{asymptotic basis of order $h$} if $hA$ 
contains every sufficiently large integer.
For example, the set of squares  is an asymptotic basis of order 4 but not of order 3.   
The set of nonnegative cubes is an asymptotic basis of order at most 7 (Linnik's theorem), and, by considering congruences modulo 9, 
an asymptotic basis of order at least 4.  
The Goldbach conjecture implies that the set of primes is an asymptotic basis of order 3.  

The modern theory of additive number theory begins with the work of Lev Genrikhovich Shnirel'man  (1905-1938).  
In an extraordinary paper~\cite{shni30}, published in Russian in 1930 and republished, 
in an expanded form~\cite{shni33}, in German in 1933, he proved that every sufficiently 
large integer is the sum of a bounded number of primes.  
Not only did Shnirel'man apply the Brun sieve, which Erd\H os subsequently  
developed into one of the most powerful tools in number theory, but he introduced 
a new density for a set of integers that is exactly the right density for the 
investigation of additive bases.
(For a survey of the classical bases in additive number theory, 
see Nathanson~\cite{nath96aa}.)

\section{Shnirel'man density and essential components}

The \emph{counting function} $A(x)$ of a set $A$ of nonnegative integers counts the number of positive integers 
in $A$ that do not exceed $x$, that is, 
\[
A(x) = \sum_{\substack{a \in A \\ 1 \leq a \leq x}} 1.
\]
The \emph{Shnirel'man density} of $A$ is
\[
\sigma(A) = \inf_{n = 1,2,\ldots} \frac{A(n)}{n}.
\]   
The sum of the sets $A$ and $B$ is the set
$A+B = \{ a+b : a \in A \text{ and } b \in B \}.$
Shnirel'man proved the fundamental  sumset inequality:
\[
\sigma(A+B) \geq \sigma(A) + \sigma(B) - \sigma(A)\sigma(B).
\]
This implies that if $\sigma(A)>0$, then $A$ is a basis of order $h$ for some $h$.
This does not apply directly to the sets of $k$th powers and the set of primes,  
which have Shnirel'man density 0.  However, it is straightforward that 
if $\sigma(A) = 0$ but $\sigma(h'A) > 0$ for some $h'$, 
then $A$ is a basis of order $h$ for some $h$.

Landau conjectured the following strengthening of  Shnirel'man's addition theorem, 
which was proved by  Mann~\cite{mann42} in 1942:
\[
\sigma(A+B) \geq \sigma(A) + \sigma(B).
\]
Artin and Scherk~\cite{arti-sche43} published a variant of Mann's proof, 
and Dyson~\cite{dyso45}, while an undergraduate at Cambridge, 
generalized Mann's inequality to $h$-fold sums.  
Nathanson~\cite{nath90c} and Heged{\" u}s, Piroska, and Ruzsa~\cite{hege-piro-ruzs98} 
have constructed examples to show that the Shnirel'man density theorems 
of Mann and Dyson are best possible.

We define the \emph{lower asymptotic density} of a set $A$ 
of nonnegative integers as follows:
\[
d_L(A) = \liminf_{n = 1,2,\ldots} \frac{A(n)}{n}.
\]   
This is a more natural density than Shnirel'man density. 
A set $A$ with asymptotic density $d_L(A) = 0$  has Shnirel'man density $\sigma(A) = 0$, 
but not conversely.  A set $A$ with asymptotic density $d_L(A) >  0$ is not necessarily 
an asymptotic basis of finite order, but $A$ is an asymptotic basis if $d_L(A) > 0$ 
and $\gcd(A) = 1$(cf. Nash and Nathanson~\cite{nath85a}).  

The set $B$  of nonnegative integers is called an \emph{essential component} if 
\[
\sigma(A+B) > \sigma(A)
\]
for every set $A$ such that $0 < \sigma(A) < 1$.  
Shnirel'man's inequality implies that every set of positive Shnirel'man density 
is an essential component.
Of course, there exist very sparse sets of zero asymptotic density 
that are not essential components.  
Khinchin~\cite{khin33} proved that the set  
of nonnegative squares is an essential component.
Note that the set of squares is a basis of order 4.  
Using an extremely clever elementary argument, Erd\H os~\cite{erdo36}, 
at the age of 22, 
proved the following considerable improvement:
Every additive basis  is an essential component.
Greatly impressed, Landau celebrated this result 
 in his 1937 Cambridge Tract 
\emph{\" Uber einige neuere {F}ortschritte der additiven {Z}ahlentheorie}~\cite{land37}.

Pl\" unnecke~\cite{plun57,plun60,plun69} and Ruzsa~\cite{ruzs89} have made important
contributions to the study of essential components.

\section{The Erd\H os-Tur\' an conjecture}

In another classic paper, published in 1941, P. Erd\H os and P. Tur\' an~\cite{erdo-tura41}  
investigated Sidon sets.  The set $A$ of nonnegative integers is a 
\emph{Sidon set} if every integer
has at most one representation as the sum of two elements of $A$.  
They concluded their paper as follows:

\begin{quotation}
Let $f(n)$ denote the number of representations of $n$ as $a_i+a_j$, \ldots.
If $f(n) > 0$ for $n > n_0$, then $\limsup f(n) = \infty.$  
Here we may mention that the corresponding result for $g(n)$, 
the number of representations of $n$ as $a_ia_j$, can be proved.
\end{quotation} 
The additive statement is still a mystery.  
The Erd\H os-Tur\' an conjecture, that the representation function 
of an asymptotic basis of order 2 is always unbounded, 
is a major unsolved problem in additive number theory.  

Many years later, in 1964, Erd\H os~\cite{erdo64a} published the proof 
of the multiplicative statement.    
This proof was later simplified by Ne\u set\u ril and R\" odl~\cite{nese-rodl85}, 
and generalized by Nathanson~\cite{nath87c}.

Long ago, while a graduate student, I searched for a counterexample 
to the Erd\H os-Tur\' an  conjecture.  
Such a counterexample might be extremal in several ways.  
It might be ``thin'' in the sense that it contains few elements.   
Every asymptotic basis of order $h$ has counting function $A(x) \gg x^{1/h}$.   
We call an additive basis of order $h$ \emph{thin} if $A(x) \ll x^{1/h}$.  Thin bases exist.  
The first examples were constructed in the 1930s by Raikov~\cite{raik37} 
and by St\" ohr~\cite{stoh37}, 
and Cassels~\cite{cass57,nath10a}  later produced another important class of examples.  

Alternatively, an asymptotic basis $A$ of order $h$ might be extremal in the sense that 
no proper subset of $A$ is an asymptotic basis of order $h$.  This means that removing 
\emph{any} element of $A$ destroys every representation of infinitely many integers.  
It is not obvious that minimal asymptotic bases exist, but I was able to construct 
asymptotic bases of order 2 that were both thin and minimal.  
Of course, none was a counterexample to the Erd\H os-Tur\' an 
conjecture.

St\" ohr~\cite{stoh55} gave the first definition of minimal asymptotic basis, 
and H\" artter~\cite{hart64} gave a non-constructive proof 
that there exist uncountably many minimal asymptotic bases 
of order $h$ for every $h \geq 2$.  

There is a natural dual to the concept of a minimal asymptotic basis.
We call a set $A$ an \emph{asymptotic nonbasis of order $h$} if is not an asymptotic 
basis of order $h$, that is, if there are infinitely many positive integers not
contained in the sumset $hA$.  An asymptotic nonbasis of order $h$ is 
\emph{maximal} if $A\cup \{b\}$ is an asymptotic basis of order $h$ for every 
nonnegative integer $b \notin A$.  The set of even nonnegative integers 
is a trivial example of a maximal nonbasis of order $h$ for every $h \geq 2$, 
and one can construct many other examples that are unions 
of the nonnegative parts of congruence classes.
It is  difficult to construct nontrivial examples.  

I discussed this and other open problems in my first paper~\cite{nath74b} 
in additive number theory.  
I did not know Erd\H os at the time, but I mailed him a preprint of the article.  
It still amazes me that he actually read this paper sent to him out of the blue 
by a completely obscure student, and he answered with a long letter 
in which he discussed his ideas about one of the problems.   
This led to correspondence, meetings, and joint work over several decades.  

\section{Extremal properties of bases}
Here is a small sample of results on minimal bases and maximal nonbases.

Nathanson and Sark\" ozy~\cite{nath89a} proved that 
if $A$ is a minimal asymptotic basis of order $h \geq 2$, then $d_L(A) \leq 1/h$.  
The proof uses Kneser's theorem~\cite{knes53} on the asymptotic density of sumsets, 
one of the most beautiful and most forgotten theorems  in additive number theory.  
A well known special case is Kneser's theorem 
for the sum of finite subsets of a finite abelian group.

Erd\H os and Nathanson~\cite{erdo-nath88a} proved that, 
for every $h \geq 2$, there exist minimal asymptotic bases of order $h$ 
with asymptotic density $1/h$.  
Moreover, for every $\alpha \in (0,1/(2h-2))$, 
there exist minimal asymptotic bases of order $h$ with asymptotic density $\alpha$.  
In particular, for every $\alpha \in (0,1/2]$, 
there exist minimal asymptotic bases of order $2$ with asymptotic density $\alpha$.

Does every asymptotic basis $A$ of order 2 contain a minimal asymptotic basis of order 2?  
Sometimes.  Let $f(n)$ count the number of representations of $n$ 
as the sum of two elements of $A$.  
If $f(n) > c\log n$ for some $c > \left(\log (4/3) \right)^{-1}$ 
and all sufficiently large $n$,    
then $A$  contains a minimal asymptotic basis of order 2 
(Erd\H os-Nathanson~\cite{erdo-nath79b}).  
This result is almost certainly not best possible.  

Does every asymptotic basis of order 2 contain a minimal asymptotic basis of order 2?  
No.  
There exists an asymptotic basis $A$ of order $2$ with the following property:  
If $S \subseteq A$,  then $A\setminus S$ is an asymptotic basis of order 2 
if and only if $S$ is finite (Erd\H os-Nathanson~\cite{erdo-nath78a}).

There exist ``trivial'' maximal asymptotic nonbases of order $h$ consisting of unions of 
arithmetic progressions~\cite{nath74b}.   
However, for every $h \geq 2$, there also exist nontrivial 
maximal asymptotic nonbases of order $h$ (Erd\H os-Nathanson~\cite{erdo-nath75a,nath77c} and Deshouillers and Grekos~\cite{desh-grek79}).

Is every asymptotic nonbasis of order $h$ a subset of a maximal asymptotic 
nonbasis of order $h$?
Sometimes.  
If $A\cup S$ is an  asymptotic nonbasis  of order $2$ 
for every finite set $S \subseteq \mathbf{N}\setminus A$, 
then $A$ is contained a maximal asymptotic nonbasis of order 2.

Is every asymptotic basis of order $h$ a subset of a maximal asymptotic nonbasis of order $h$?
No.
Hennefeld~\cite{henn77} proved that, for every $h \geq 2$ there exists 
an asymptotic nonbasis $A$ of order $h$ such that, 
if $S \subseteq \mathbf{N}\setminus A$, then $A\cup S$ is an  asymptotic 
nonbasis $A$ of order $h$ if and only if the set $\mathbf{N}\setminus (A \cup S)$ is infinite.

Investigating extremal properties of additive bases is like exploring for 
new plant species in the Amazon rain forest.  Much has been collected, 
but much more is unimagined.  
The following results about oscillations of bases and nonbases 
appear in~\cite{erdo-nath75b,erdo-nath76}.

There exists a minimal asymptotic basis of order 2 such that 
$A\setminus \{x\}$ is a maximal asymptotic nonbasis of order 2 for every $x \in A$.

There exists a maximal asymptotic nonbasis of order 2 such that 
$A\cup \{ y \}$ is a minimal asymptotic basis of order 2 
for every $y \in \mathbf{N} \setminus A$.

There exists a partition of the nonnegative integers into disjoint sets $A$ and $B$ such 
that $A$ is a minimal asymptotic basis of order 2 and $B$ is a maximal 
asymptotic nonbasis of order 2.  

There exists a partition of the nonnegative integers into disjoint sets $A$ and $B$ 
that oscillate in phase from minimal asymptotic basis of order 2 to  maximal 
asymptotic nonbasis of order 2 as random elements are moved 
from the basis to the nonbasis, infinitely often.  

It is an open problem to extend these results to asymptotic bases of order $h \geq 3$.  
For a survey of extemal problems in additive number theory, see Nathanson~\cite{nath15aa}.

Bases fascinated Erd\H os, but the subject is still obscure, an obsession of the few.
I wrote my first paper~\cite{nath74b} on minimal bases and maximal nonbases 
in additive number theory in the summer of 1970, while a graduate student, 
visiting the Weizmann Institute in Israel.   
I was thinking about the Erd\H os -Tur\' an conjecture, and studied the 
foundational book  \emph{Sequences}~\cite{halb-roth66} by  Halberstam and 
Roth  (now, unfortunately, out of print).    
Returning to Weizmann to lecture in 2001, I looked for the book in the library.  
In the ancient, pre-computer era, a library book had a 3 $\times$ 5 card.  
When you checked out a book, you wrote your name on the card.  
You could see who had read any book.  
After 30 years, mine was still the only name on the card for \emph{Sequences}.

\def\cprime{$'$} \def\cprime{$'$} \def\cprime{$'$} \def\cprime{$'$}
\providecommand{\bysame}{\leavevmode\hbox to3em{\hrulefill}\thinspace}
\providecommand{\MR}{\relax\ifhmode\unskip\space\fi MR }
\providecommand{\MRhref}[2]{%
  \href{http://www.ams.org/mathscinet-getitem?mr=#1}{#2}
}
\providecommand{\href}[2]{#2}


\begin{thebibliography}{10}

\bibitem{arti-sche43}
E.~Artin and P.~Scherk, \emph{On the sum of two sets of integers}, Annals Math.
  \textbf{44} (1943), 138--142.

\bibitem{cass57}
J.~W.~S. Cassels, \emph{{\" U}ber {B}asen der nat{\" u}rlichen {Z}ahlenreihe},
  Abhandlungen Math. Seminar Univ. Hamburg \textbf{21} (1975), 247--257.

\bibitem{desh-grek79}
J.-M. Deshouillers and G.~Grekos, \emph{Propri{\' e}t{\' e}s extr{\' e}ales de
  bases additives}, Bull. Soc. Math. France \textbf{107} (1979), 319--335.

\bibitem{dyso45}
F.~J. Dyson, \emph{A theorem on the densities of sets of integers}, J. London
  Math. Soc. \textbf{20} (1945), 8--14.

\bibitem{erdo-tura41}
P.~Erd\H{o}s and P.~Tur\'an, \emph{On a problem of {S}idon in additive number
  theory, and on some related problems}, J. London Math. Soc. \textbf{16}
  (1941), 212--215.

\bibitem{erdo36}
P.~Erd{\H o}s, \emph{On the arithmetical density of the sum of two sequences,
  one of which forms a basis for the integers}, Acta Arith. \textbf{1} (1936),
  197--200.

\bibitem{erdo64a}
\bysame, \emph{On the multiplicative representation of integers}, Israel J.
  Math. \textbf{2} (1964), 251--261.

\bibitem{erdo-nath75a}
P.~Erd{\H o}s and M.~B. Nathanson, \emph{Maximal asymptotic nonbases}, Proc.
  Amer. Math. Soc. \textbf{48} (1975), 57--60.

\bibitem{erdo-nath75b}
\bysame, \emph{Oscillations of bases for the natural numbers}, Proc. Amer.
  Math. Soc. \textbf{53} (1975), 253--358.

\bibitem{erdo-nath76}
\bysame, \emph{Partitions of the natural numbers into infinitely oscillating
  bases and nonbases}, Commentarii Math. Helvet. \textbf{52} (1976), 171--182.

\bibitem{nath77c}
P.~Erd{\H{o}}s and M.~B. Nathanson, \emph{Nonbases of density zero not
  contained in maximal nonbases}, J. London Math. Soc. (2) \textbf{15} (1977),
  no.~3, 403--405.

\bibitem{erdo-nath78a}
P.~Erd{\H o}s and M.~B. Nathanson, \emph{Sets of natural numbers with no
  minimal asymptotic bases}, Proc. Amer. Math. Soc. \textbf{70} (1978),
  100--102.

\bibitem{erdo-nath79b}
\bysame, \emph{Systems of distinct representatives and minimal bases in
  additive number theory}, Number Theory, Carbondale, 1979 (Heidelberg) (M.~B.
  Nathanson, ed.), Lecture Notes in Mathematics, vol. 751, Springer-Verlag,
  1979, pp.~89--107.

\bibitem{erdo-nath88a}
\bysame, \emph{Minimal asymptotic bases with prescribed densities}, Illinois J.
  Math. \textbf{32} (1988), 562--574.

\bibitem{halb-roth66}
H.~Halberstam and K.~F. Roth, \emph{{Sequences, Vol. 1}}, Oxford University
  Press, Oxford, 1966, Reprinted by Springer-Verlag, Heidelberg, in 1983.

\bibitem{hart64}
E.~H{\"a}rtter, \emph{Eine {B}emerkung \"uber periodische {M}inimalbasen f\"ur
  die {M}enge der nichtnegativen ganzen {Z}ahlen}, J. Reine Angew. Math.
  \textbf{214/215} (1964), 395--398.

\bibitem{hege-piro-ruzs98}
P.~Heged{\" u}s, G.~Piroska, and I.~Z. Ruzsa, \emph{On the {S}chnirelmann
  density of sumsets}, Publ. Math. Debrecen \textbf{53} (1998), no.~3--4,
  333--345.

\bibitem{henn77}
J.~Hennefeld, \emph{Asymptotic nonbases which are not subsets of maximal
  aymptotic nonbases}, Proc. Amer. Math. Soc. \textbf{62} (1977), 23--24.

\bibitem{khin33}
A.~Ya. Khinchin, \emph{{\"U}ber ein metrisches {P}roblem der additiven
  {Z}ahlentheorie}, Mat. Sbornik N.S. \textbf{10} (1933), 180--189.

\bibitem{knes53}
M.~Kneser, \emph{Absch{\"a}tzungen der asymptotischen {D}ichte von
  {S}ummenmengen}, Math. Z. \textbf{58} (1953), 459--484.

\bibitem{land37}
E.~Landau, \emph{{\"U}ber einige neuere {F}ortschritte der additiven
  {Z}ahlentheorie}, Cambridge University Press, Cambridge, 1937.

\bibitem{mann42}
H.~B. Mann, \emph{A proof of the fundamental theorem on the density of sums of
  sets of positive integers}, Annals Math. \textbf{43} (1942), 523--527.

\bibitem{nath85a}
J.~C.~M. Nash and M.~B. Nathanson, \emph{Cofinite subsets of asymptotic bases
  for the positive integers}, J. Number Theory \textbf{20} (1985), no.~3,
  363--372.

\bibitem{nath74b}
M.~B. Nathanson, \emph{Minimal bases and maximal nonbases in additive number
  theory}, J. Number Theory \textbf{6} (1974), 324--333.

\bibitem{nath87c}
\bysame, \emph{Multiplicative representations of integers}, Israel J. Math.
  \textbf{57} (1987), no.~2, 129--136.

\bibitem{nath90c}
\bysame, \emph{Best possible results on the density of sumsets}, Analytic
  number theory (Allerton Park, IL, 1989), Progr. Math., vol.~85, Birkh\"auser
  Boston, Boston, MA, 1990, pp.~395--403.

\bibitem{nath96aa}
\bysame, \emph{{Additive Number Theory: The Classical Bases}}, Graduate Texts
  in Mathematics, vol. 164, Springer-Verlag, New York, 1996.

\bibitem{nath10a}
M.~B. Nathanson, \emph{Cassels bases}, {Additive Number Theory}, Springer, New
  York, 2010, pp.~259--285.

\bibitem{nath15aa}
M.~B. Nathanson, \emph{{Additive Number Theory: Extremal Problems and the
  Combinatorics of Sumsets}}, Graduate Texts in Mathematics, Springer, New
  York, to appear.

\bibitem{nath89a}
M.~B. Nathanson and Andr{\'a}s S{\'a}rk{\"o}zy, \emph{On the maximum density of
  minimal asymptotic bases}, Proc. Amer. Math. Soc. \textbf{105} (1989),
  31--33.

\bibitem{nese-rodl85}
J.~Ne{\u s}et{\u r}il and V.~R{\" o}dl, \emph{Two proofs in combinatorial
  number theory}, Proc. Amer. Math. Soc. \textbf{93} (1985), 185--188.

\bibitem{plun57}
H.~Pl{\"u}nnecke, \emph{{\"U}ber ein metrisches {P}roblem der additiven
  {Z}ahlentheorie}, J. reine angew. Math. \textbf{197} (1957), 97--103.

\bibitem{plun60}
\bysame, \emph{{\"U}ber die {D}ichte der {S}umme zweier {M}engen, deren eine
  die dichte null hat}, J. reine angew. Math. \textbf{205} (1960), 1--20.

\bibitem{plun69}
\bysame, \emph{Eigenschaften und absch{\"a}tzungen von wirkungsfunktionen},
  vol.~22, Berichte der Gesellschaft f{\"u}r Mathematik und Datenverarbeitung,
  Bonn, 1969.

\bibitem{raik37}
D.~Raikov, \emph{{\" U}ber die {B}asen der nat{\" u}rlichen {Z}ahlentreihe},
  Mat. Sbornik N. S. 2 \textbf{44} (1937), 595--597.

\bibitem{ruzs89}
I.~Z. Ruzsa, \emph{An application of graph theory to additive number theory},
  Scientia, Ser. A \textbf{3} (1989), 97--109.

\bibitem{shni30}
L.~G. Shnirel'man, \emph{On the additive properties of integers}, Izv. Donskovo
  Politekh. Inst. Novocherkasske \textbf{14} (1930), 3--27.

\bibitem{shni33}
\bysame, \emph{{\"U}ber additive {E}igenschaften von {Z}ahlen}, Math. Annalen
  \textbf{107} (1933), 649--690.

\bibitem{stoh37}
A.~St{\"o}hr, \emph{Eine {B}asis {$h$-O}rdnung f{\" u}r die {M}enge aller
  nat{\" u}rlichen {Z}ahlen}, Math. Zeit. \textbf{42} (1937), 739--743.

\bibitem{stoh55}
A.~St{\"o}hr, \emph{Gel\"oste und ungel\"oste {F}ragen \"uber {B}asen der
  nat\"urlichen {Z}ahlenreihe. {I}, {II}}, J. Reine Angew. Math. \textbf{194}
  (1955), 40--65, 111--140.

\end{thebibliography}
\end{document}